\documentclass[11pt,reqno,oneside]{amsart}
\usepackage{amsmath,amssymb,cite,epsfig,xspace,alltt,verbatim}

\numberwithin{equation}{section}

\theoremstyle{plain}
\newtheorem{theorem}{Theorem}[section]
\newtheorem{lemma}{Lemma}[section]

\theoremstyle{definition}

\newtheorem{problem}{Problem}[section]

\theoremstyle{remark}

\newcommand{\Real}{\mathbb R}

\newcommand{\Int}{\mathbb Z}
\newcommand{\LP}{\mathcal L}
\newcommand{\I}{\mathcal I}
\newcommand{\K}{\kappa}
\newcommand{\psib}{\mathbf \Psi}
\newcommand{\phib}{\mathbf \Phi}

\newcommand{\eb}{\mathbf e}
\newcommand{\fb}{\mathbf f}

\newcommand{\rb}{\mathbf z}
\newcommand{\sbf}{\mathbf s}
\newcommand{\yb}{\mathbf y}
\newcommand{\zb}{\mathbf z}

\newcommand{\ud}{\mathrm{d}}
\newcommand{\A}{{\mathcal A}}
\newcommand{\E}{{\mathcal E}}
\newcommand{\W}{{\mathcal W}}
\newcommand{\Po}{{\mathcal E}_{ext}}
\newcommand{\const}{CQC}
\newcommand{\define}{:=}
\newcommand{\nr}{\nabla_{\negthickspace \rb}}
\def\maxnorm#1{\left|\left|#1\right|\right|_\infty}
\newcommand{\half}{\frac{1}{2}}
\DeclareMathOperator*{\argmin}{argmin}

\newcommand{\pd}[2]{\frac{\partial #1}{\partial #2}}

\begin{document}

%%%%%%%%%%%%%
% Topmatter %
%%%%%%%%%%%%%
\title[Iterative Solution of the QC Equilibrium Equations with Continuation]
{Iterative Solution of the Quasicontinuum Equilibrium
Equations with Continuation}
%%%%%%%%%%%
% Authors %
%%%%%%%%%%%
\author{Matthew Dobson}
\author{Mitchell Luskin}

%%%%%%%%%%%%%
% Addresses %
%%%%%%%%%%%%%
\address{Matthew Dobson\\
School of Mathematics \\
University of Minnesota \\
206 Church Street SE \\
Minneapolis, MN 55455 \\
U.S.A.}
\email{dobson@math.umn.edu}

\address{Mitchell Luskin \\
School of Mathematics \\
University of Minnesota \\
206 Church Street SE \\
Minneapolis, MN 55455 \\
U.S.A.}
\email{luskin@umn.edu}

%%%%%%%%%%%%%%%%%%
% Acknowledgments %
%%%%%%%%%%%%%%%%%%
\thanks{This work was supported in part by
 DMS-0304326
  and by the University of Minnesota Supercomputing Institute.
  This work is also based on
work supported by the Department of Energy under Award Number
DE-FG02-05ER25706.
}

%%%%%%%%%%%%
% Keywords %
%%%%%%%%%%%%

%%%%%%%%%%%%%%%%%
% Subject class %
%%%%%%%%%%%%%%%%%

%%%%%%%%
% Date %
%%%%%%%%
\date{\today}

%%%%%%%%%%%%
% Abstract %
%%%%%%%%%%%%
\begin{abstract}
We give an analysis of a continuation algorithm for the numerical solution of
the force-based quasicontinuum equations.  The approximate solution of the
force-based quasicontinuum equations is computed by an iterative method using
an energy-based quasicontinuum approximation as the preconditioner.

The analysis presented in this paper is used to determine an efficient strategy
for the parameter step size and number of iterations at each parameter value to
achieve a solution to a required tolerance.  We present computational results
for the deformation of a Lennard-Jones chain under tension to demonstrate the
necessity of carefully applying continuation to ensure that the computed
solution remains in the domain of convergence of the iterative method as the
parameter is increased.  These results exhibit fracture before the actual load
limit if the parameter step size is too large.
\end{abstract}

%%%%%%%%
% Body %
%%%%%%%%

\maketitle
\section{Introduction}
Quasicontinuum (QC) approximations reduce the computational
complexity of a material simulation by reducing the degrees of
freedom used to describe a configuration of atoms and by giving
approximate equilibrium equations on the reduced degrees of
freedom~\cite{pinglin03, pinglin05, legollqc05, ortnersuli,
OrtnerSueli:2006d, minge07, tadmor_qc_first, knaportiz,e06, e05, mill02,
rodney_gf, miller_indent, curtin_miller_coupling,
jacobsen04,dobs07}.  For crystalline materials, there are
typically a few small regions with highly non-uniform structure
caused by defects in the material which are surrounded by large
regions where the local environment of atoms varies slowly.
The idea of QC is to replace these slowly varying regions with a
continuum model and couple it directly to the atomistic model
surrounding the defects.  The material's position is described by a set
of representative atoms that are in one-to-one correspondence
with the lattice atoms in the atomistic regions but reduce the
degrees of freedom in the
continuum regions.

Quasi-static computations in material simulations explore
mechanical response under slow external loading by fully relaxing
the material at each step of a parameterized path of external conditions.
Such simulations can model nano-indentation, stress-induced phase
transformations, and many other material processes. The characteristic feature
when using this technique is that the process to be modeled occurs slowly
enough that dynamics are assumed to play no role in determining the relaxed
state.
This paper focuses on applying continuation
techniques~\cite{bank82,kell77,doed97} to the
nonlinear equilibrium equations of the
force-based quasicontinuum approximation (QCF).

\subsection{Choosing a Quasicontinuum Approximation}
There are many choices available for the interaction among the
representative atoms, especially between those in the atomistic
and continuum regions, which has led to the development of a
variety of quasicontinuum approximations.  Criteria for
determining a good choice of approximation for a given problem
are still being developed.  Algorithmic simplicity and
efficiency is certainly important for implementation and
application, but concerns about accuracy have led to the search
for consistent schemes.  Since the forces on all of the atoms in a
uniformly strained lattice are zero,
we define a quasicontinuum approximation to be consistent if
there are no forces on the representative atoms for a lattice
that has been deformed by a uniform strain.  We note that the
atomistic to continuum interface is typically where consistency
fails for a QC approximation~\cite{dobs07,mill02,jacobsen04,e06}
as is the case for the original QC
method~\cite{tadmor_qc_first}.

For static problems, QCF is an attractive choice for
quasicontinuum approximation because it is a consistent
scheme~\cite{dobs07,mill02} that is algorithmically simple:
the force on each representative atom is given by either
an atomistic calculation or a continuum finite element calculation.
The algorithmic simplicity of the force-based
quasicontinuum method has allowed it to be implemented with
adaptive mesh refinement and atomistic to continuum model selection
algorithms
~\cite{mill02,curtin_miller_coupling,arndtluskin07a,arndtluskin07b,ArndtLuskin:2007c,PrudhommeBaumanOden:2005,OdenPrudhommeRomkesBauman:2005}.
The trade-off for the consistency and algorithmic simplicity
of QCF is that it does not give a conservative force field,
although it is close to a conservative field~\cite{dobs07}.
Thus, QCF is a method to approximate forces, rather than a
method to approximate the energy.

Quasicontinuum energies have been proposed that
utilize special energies for atoms in an interfacial region~\cite{jacobsen04,e06},
and corresponding conditions for consistency have been satisfied for planar interfaces~\cite{e06}.
However, there is
currently no known consistent quasicontinuum energy that allows general
nonplanar atomistic to continuum interfaces and mesh coarsening in the continuum region
(other than the computationally intensive
constrained quasicontinuum energy discussed in Section~\ref{secconstrain}).
We will, however, investigate the use of quasicontinuum energies
as preconditioners for the force-based quasicontinuum approximation.

\subsection{Solving Equilibrium Equations by Continuation}
Our goal is to efficiently approximate the
solution $\zb(s)$ to the QCF equilibrium
equations
\begin{equation*}
F^{QCF}(\zb(s),s) + f(s) = 0,\qquad s\in[0,1],
\end{equation*}
where $\zb(s) \in \Real^n$ are the coordinates of
the representative atoms that describe the material and the
parameter $s \in [0,1]$ represents the change in external
loads such as an indenter position or an applied force.
Using continuation,
we start from $\zb(0)$ which is usually easy to find (such as
a resting position) and follow the solution branch by incrementing $s$ and
looking for a solution $\zb(s)$ in a neighborhood of the previous solution.
The continuation approach
that we analyze in this paper has been used to obtain
computational solutions to materials deformation
problems~\cite{mill02,curtin_miller_coupling} and is
implemented in the multidimensional QC code~\cite{qcmethod}.

The approach that we will investigate in this
paper uses constant extrapolation in $s$ to obtain initial states
for the iterative solution of $F^{QCF}(\zb,s)$
at a sequence of load steps $s_q$ where
$0=s_0 \leq s_1 \leq \dots \leq s_Q=1$ and
solves the iterative equations
using a preconditioner force $F^{QCE}$ that comes from a quasicontinuum energy
$\E^{QCE}(\zb,s),$
that is, $F^{QCE}(\zb, s) = -\frac{\partial
\E^{QCE}(\zb,s)}{\partial \zb}.$
We will focus our analysis on a specific preconditioner later,
but the splitting and subsequent analysis
works for any choice of quasicontinuum energy.
The outer iteration at a fixed step $s_q$ is
given by
\begin{equation}
\label{precona}
\begin{split}
F^{QCE}(\zb_q^{p+1},s_{q})&=
F^{QCE}(\zb_q^{p},s_{q})-F^{R}(\zb_q^{p},s_{q})\\
&=-f(s_{q})-F^{G}(\zb_q^{p},s_{q}),
\qquad p=0,\dots,P_q-1,\\
\zb_q^0&=\zb_{q-1},
\end{split}
\end{equation}
where
\[
F^{R}(\zb,s)\define
f(s)+F^{QCF}(\zb,s)
\]
is the residual force and
\[
F^{G}(\zb,s)\define
F^{QCF}(\zb,s)-F^{QCE}(\zb,s)
\]
is sometimes called a ``ghost force correction''
in the mechanics literature~\cite{mill02}.  We will consider
preconditioner forces  $F^{QCE}(\zb,s)$ that differ from
$F^{QCF}(\zb,s)$ only in atomistic to continuum interfacial regions
so that the ghost force correction is inexpensive to compute.  Since
the preconditioner forces come from an energy, the outer
iteration equations \eqref{precona} for $\zb_q^{p+1}$ can
be solved by an inner iteration that finds the local minimum
in a neighborhood of $\zb_q^p$ for
\begin{equation}\label{inner}
\left[\E^{QCE}(\zb,s_{q})
-\left(f(s_q)+
F^{G}(\zb_q^{p},s_{q})\right)\cdot\zb\right]
\end{equation}
using an energy minimization method starting with initial guess $\zb =
\zb_q^p.$

We give an analysis to optimize the computational efficiency of
the continuation algorithm~\eqref{precona} by varying the parameter step
size, $h_q=s_q-s_{q-1},$ and the number of
outer iterations, $P_q.$  Our analysis first considers the goal
of computing an approximation of $\zb(s)$ uniformly for
$s\in[0,1]$ to a given tolerance, $\epsilon.$  The proposed
strategy selects the step size $h_q$ so that the initial
iterates $\zb_q^{0}$ are in the domain of convergence of the
outer iteration and so that the tolerance is achieved by the
continuous, piecewise linear interpolant of the solution at each
parameter $s_q.$ The required accuracy at the $s_q$ is achieved
by the fast convergence of the iteration~\eqref{precona}.  As
$\epsilon\to 0,$ our analysis gives that $h_q\to 0$ and
$P_q\to\infty$ for all $q=1,\dots,Q,$ so that an efficient way
to achieve increased accuracy uses a balance between small step
size for accurate interpolation and a large number of iterations
per step.

We then consider the goal of efficiently computing the
final state $\zb(1)$ to a given tolerance, $\epsilon.$
For this goal, the result of our analysis states that
an efficient strategy fixes the number of outer iterations to $P_q=1$
at all but the final step
and takes the largest possible steps, $h_q,$
such that the initial guesses, $\zb_q^{0},$
remain in the domain of convergence of the iteration.
This strategy determines the number of steps, $Q,$
independently of $\epsilon.$
The required tolerance, $\epsilon,$ is then achieved at
$s=1$ by doing sufficiently many iterations, $P_Q>1.$
In this case, only $P_Q\to\infty$ as $\epsilon\to 0.$

We give numerical results for the deformation of a
Lennard-Jones chain under tension that demonstrate
the importance of selecting the
step size, $h_q,$
and number of outer iterations, $P_q,$ so that the
iterates $\zb_q^{p}$
remain in the domain of convergence of the iteration.
The numerical experiment
shows that our algorithm diverges (the chain spuriously undergoes fracture)
if we attempt to solve for the deformation corresponding
to a load near the limit load by a step size $h_1=1.$

We give a derivation of the force-based quasicontinuum
approximation and the energy-based preconditioner in
Section~\ref{sec:qc}.  In Section~\ref{sec:iter}, we analyze the
equilibrium equations and their iterative solution.  In
Section~\ref{tension}, we apply Theorem~\ref{thm:contract} to a
Lennard-Jones chain under tension to obtain bounds on the
initial state that guarantee convergence of our iterative method
to the equilibrium state, to obtain convergence results for our
iterative method, and to demonstrate the need for continuation
by the showing that the chain can undergo fracture if we begin
our iteration outside the prescribed neighborhood.
Section~\ref{sec:cont} presents the continuation method and
Sections~\ref{sec:uniform} and \ref{sec:contend} give an
analysis to guide the development of an efficient algorithm
using the quasicontinuum iteration.  We collect the proofs of
several lemmas in a concluding Appendix~\ref{appendix}.

\section{Quasicontinuum Approximations}
\label{sec:qc}

This section describes a model for a one-dimensional chain of
atoms and a sequence of approximations that lead to the
force-based and energy-based quasicontinuum approximations.
While a one-dimensional model is very limited in the type of
defects it can exhibit, its study illustrates many of the
theoretical and computational issues  of QC approximations.

We treat the case where atomistic interactions are governed by a
pairwise classical potential $\phi(r),$ where $\phi$ is defined
for all $r>0.$ A short-range cutoff for the potential is a good
approximation for many crystals, and in the following analysis
we use a second-neighbor (next-nearest neighbor) cutoff, as this
gives the simplest case in which the atomistic and continuum
models are distinct~\cite{dobs07}.

\subsection{The Fully Atomistic Model}
Denote the positions of the atoms in a linear chain by $y_i$ for
$i=-M,\dots,M+1$ where $y_i < y_{i+1}$ and denote the position
where the right-hand end of the chain is fixed by
$y_{M+1} = \hat y_{M+1}(s)$ for a parameter $s\in[0,1].$ The
second-neighbor energy for the chain is
then given by
\begin{equation}
\label{atomTot}
\E^a(\yb,s) \define \sum_{i=-M}^M \phi(y_{i+1} - y_i)
    + \sum_{i=-M}^{M-1} \phi(y_{i+2} - y_i)
\end{equation}
where $\yb\define(y_{-M},\dots,y_{M+1}).$ We also assume that the
chain is subject to an external potential energy which we assume
for simplicity to have the form
\begin{equation*}
\Po^a(\mathbf y,s)\define -\sum_{i=-M}^{M}\tilde{f}_i(s) y_i.
\end{equation*}
Section \ref{tension} describes a numerical
example with a boundary dead-load given by
$\tilde{f}_{-M}(s)\ne 0$ and
$\tilde{f}_i(s)=0$ for all interior atoms $i=-M+1,\dots,M.$

We want to find local minima of the total energy,
\begin{equation}
\label{atomSol}
\E^{a}_{total}\define\E^a(\yb,s)+\Po^a(\mathbf y,s)
\end{equation}
subject to the boundary constraint $y_{M+1} = \hat y_{M+1}(s).$
The equilibrium equation for the fully atomistic system
\eqref{atomSol} is given by
\begin{equation*}
F_i^{a}(\yb(s),s) + \tilde f_i(s)=0,\qquad i=-M,\dots,M,
\end{equation*}
where the atomistic force is given by
\begin{equation*} 
\begin{split}
F_{i}^a(\mathbf y,s) \define -\frac{\partial \E^a(\yb,s)}{\partial y_i}
&= \left[\phi'(y_{i+1} - y_{i}) + \phi'(y_{i+2} - y_i)\right] \\
&\qquad - \left[\phi'(y_i - y_{i-1}) + \phi'(y_{i}-y_{i-2})\right],
\end{split}
\end{equation*}
for $i=-M,\dots,M$ where the
terms $\phi'(y_i -y_j)$
above and in the following should be understood to be zero
for $i\notin \{-M,\dots,M+1\}$ or $j\notin \{-M,\dots,M+1\}.$
In the remainder of this section, we will not explicitly
denote the dependence on the parameter $s.$

\subsection{The Constrained Quasicontinuum Approximation}
\label{secconstrain}
The constrained quasicontinuum approximation finds approximate
minimum energy configurations of \eqref{atomSol} by selecting
 a subset of the atoms to act as representative
atoms and interpolating the remaining atom positions via
piecewise linear interpolants in the reference configuration.
We denote by $a_0$ the ground state lattice constant
for the potential $\phi(r)$ with a second-neighbor cutoff, that
is,
\begin{equation*}
a_0 \define \argmin \phi(r)+\phi(2r)
\end{equation*}
(see Section~\ref{sec:iter} or \cite{dobs07}).
We then set the reference positions
of the atoms in the chain to be
\begin{equation*}
x_i \define i a_0 \quad \textrm{for } i=-M,\dots,M+1.
\end{equation*}

We let $z_j\define y_{\ell_j}$ denote the representative
atom positions where $j = -N,\dots,N+1,$
and where $\ell_{-N}=-M,$ $\ell_{N+1}=M+1,$ and $\ell_j<\ell_{j+1}.$
We are interested in developing methods
for $N \ll M.$ We can obtain the positions of all atoms
$y_i$ from the positions of representative atoms $z_j$ by
\begin{equation*}
y_i(\zb)=\sum_{j=-N}^{N+1}S_j(x_i)z_j\quad\text{for } i=-M,\dots,M+1,
\end{equation*}
where $\zb\define(z_{-N},\dots,z_{N+1})$ and the $S_j(x)$ are the continuous,
piecewise linear ``shape'' functions for the mesh constructed
from the reference coordinates $x_{\ell_j}$
of the representative atoms, more precisely,
\begin{equation} \label{eq:basisfcn}
  S_j(x)\define \begin{cases}
  0,                   &\text{if }x\leq x_{\ell_{j-1}},\\
  (x-x_{\ell_{j-1}} )/(x_{\ell_{j}} - x_{\ell_{j-1}}),
    &\text{if }x_{\ell_{j-1}}
  <x\leq x_{\ell_{j}},\\
  (x_{\ell_{j+1}} -x)/(x_{\ell_{j+1}} - x_{\ell_{j}}),
  &\text{if }x_{\ell_{j}}<x\leq x_{\ell_{j+1}},\\
  0,                   &\text{if }x>x_{\ell_{j+1}}.
\end{cases}
\end{equation}

The {\em constrained quasicontinuum energy} is then given by
\begin{equation*}
\E^{\const}(\zb)\define \E^a(\yb(\zb)),
\end{equation*}
and the {\em constrained external potential energy} is given by
\begin{equation*}
\Po^{\const}(\zb)\define \Po^a(\yb(\zb)).
\end{equation*}
Using \eqref{eq:basisfcn} and
the chain rule, we obtain the conjugate atomistic force, that is, the force
on the reduced degrees of freedom induced by the atomistic
forces.   We find that
\begin{equation*}
\begin{split}
F^{\const}_j(\zb) &\define -\frac{\partial
\E^{\const}(\zb)}{\partial z_j} \\
&=
\sum_{i=0}^{\nu_{j-1}}
 \left(\frac{\nu_{j-1} -i}{\nu_{j-1}}\right)F^{a}_{\ell_j-i}(\yb(\zb))
 +\sum_{i=1}^{\nu_{j}} \left(\frac{\nu_{j} - i}{\nu_{j}}\right)
F^{a}_{\ell_j+i}(\yb(\zb))
\end{split}
\end{equation*}
for $j=-N,\dots,N,$ and the conjugate external force is given by
\begin{equation}\label{quasifa}
f_{j}\define -\frac{\partial \Po^{\const}(\zb)}{\partial z_j}=\sum_{i=0}^{\nu_{j-1}}
 \left(\frac{\nu_{j-1} -i}{\nu_{j-1}}\right)\tilde{f}_{\ell_j-i}
 +\sum_{i=1}^{\nu_{j}} \left(\frac{\nu_{j} - i}{\nu_{j}}\right)
 \tilde{f}_{\ell_j+i},
\end{equation}
for $j=-N,\dots,N,$ where
\begin{equation*}
\nu_j \define \ell_{j+1} - \ell_{j}
\end{equation*} is the number
of atoms between $z_{j}$
and $z_{j+1}$ (the end atoms are only counted half).
The equilibrium equations for the total constrained quasicontinuum energy,
\[
\E^{\const}_{total}(\zb) \define \E^{\const}(\zb)+\Po^{\const}(\zb)
\]
are then given by
\begin{equation*}
F_j^{\const}(\zb) + f_j=0,\qquad j=-N,\dots,N.
\end{equation*}

The constrained quasicontinuum approximation is attractive since it gives
conservative forces and since it is the only known quasicontinuum energy that
is consistent when generalized to multidimensional approximations~\cite{e06}.
The constrained quasicontinuum approximation is also
attractive since its conjugate forces~\eqref{quasifa} are located at only $2N$
representative atoms; however, we must still compute the forces at all $2M$
atoms which makes it computationally infeasible.  Some computational savings
can be made in the interior of large elements by separating the energy
computations into element energy plus surface energy; however, in higher
dimensions the large number of atoms near element boundaries makes the
constrained quasicontinuum approximation impractical.
Finally, it is attractive because its approximation error comes only
from the restriction to linear deformations within the element making
it possible to analyze using classical finite element error analysis.

\subsection{The Local Quasicontinuum Energy}
We now recast the constrained approximation in terms of continuum mechanics
to introduce the {\em local quasicontinuum energy} which is
simply a continuous, piecewise linear approximation of a hyperelastic continuum model where the
strain-energy density is derived from the atomistic potential,
$\phi(r).$ This energy efficiently approximates the conjugate force at
the representative atoms.
We have that~\cite{dobs07}
\begin{equation*}
\begin{split}
{\mathcal E}^{\const}(\mathbf z) &=\sum_{j=-N}^{N}
  L_j W\left(D_j\right) + \mathcal S_b(D_{-N}) \\
&\qquad  +\sum_{j=-N+1}^{N}\mathcal S\left(D_{j-1}, D_j\right)+
  \mathcal S_b(D_{N}),
\end{split}
\end{equation*}
where
\begin{equation*}
L_j \define x_{\ell_{j+1}}-x_{\ell_j}\quad
\text{and}\quad D_j\define\frac{z_{j+1} - z_{j}}{x_{\ell_{j+1}}-x_{\ell_j}}
\end{equation*}
are the length and deformation gradient of the $j$th element,
and
\begin{equation*}
W(D) \define \frac{\phi(Da_0)+\phi(2Da_0)}{a_0},
\end{equation*}
is the strain-energy density for an
infinite atomistic chain with the uniform lattice spacing $D a_0.$
Here~\cite{dobs07}
\begin{equation*}
\begin{split}
\mathcal S_b(D)=&-\frac{1}{2}\phi(2Da_0),\\
\mathcal S(D_{1},D_2) =& - \frac{1}{2}\phi(2D_{1} a_0) +
 \phi(D_{1}a_0 + D_2 a_0)
    - \frac{1}{2}\phi(2D_2 a_0),\\
\end{split}
\end{equation*}
can be considered to be a surface energy
and an interfacial energy respectively.

Since $\mathcal S(D_{j-1},D_j)$ is a second divided difference,
the interfacial energy is small in regions where
the strain is slowly varying.  We obtain the {\em local quasicontinuum
energy} by neglecting the interfacial energy and surface energy
to obtain
\begin{equation*}
\E^L(\zb)\define\sum_{j=-N}^{N}
  L_j W\left(D_j\right),
  \end{equation*}
and we have the corresponding conjugate atomistic force
\begin{equation*}
F^{L}_j(\zb)\define -\frac{\partial \E^L(\zb)}{\partial z_j}
=\frac{\partial W}{\partial D}(D_j)
-\frac{\partial W}{\partial D}(D_{j-1}),\qquad
j=-N,\dots, N.
\end{equation*}

We note that $F^{L}_j(\zb)$
depends only
on $z_{j-1},$ $z_j,$ and $z_{j+1}.$  This approximation is
computationally feasible since the work to compute all the forces
is proportional to N.  The approximation error now has two components:
the linearization within each element that is inherited from the constrained quasicontinuum
approximation plus the operator error incurred by ignoring interfacial terms.
In cases where the deformation gradient $D_j$ is slowly varying on the scale of
the representative atom mesh, both sources of error will be small and
the local approximation will be highly accurate, as is expected for a sufficiently refined
finite element continuum model.  Mesh refinement can be used to reduce
both sources of error, but even mesh refinement to the atomistic scale
cannot remove the interfacial error in the neighborhood of defects since the deformation
gradient varies rapidly on the atomistic scale.  Thus,
the atomistic model must be retained near defects to obtain sufficient
accuracy.

\subsection{The Force-based Quasicontinuum Approximation}
 We can obtain a
quasicontinuum approximation that is accurate in regions
where the deformation gradient $D_j$ is rapidly varying,
such as in the neighborhood of defects, and maintains
the efficiency of the local quasicontinuum method by combining
them in the {\em force-based quasicontinuum} approximation (QCF).
In QCF, we partition the chain into ``atomistic'' and ``continuum''
representative atoms and define the force on each
representative atom to be the
force that would result if the whole approximation was of its
respective type, that is,
\begin{equation}
\label{force}
F_{j}^{QCF}(\mathbf z)\define
\begin{cases}
F_j^{a}(\zb) & \textrm{if representative atom }
j \textrm{ is atomistic}, \\
F_j^{L}(\zb)&
\textrm{if representative atom }
j \textrm{ is continuum}.
\end{cases}
\end{equation}
 With this convention, for example,
the forces on a continuum representative atom are determined solely by the
adjacent degrees of freedom regardless of how close any
atomistic representative atoms may be.

\begin{figure}[htb]
\includegraphics[width=\textwidth]{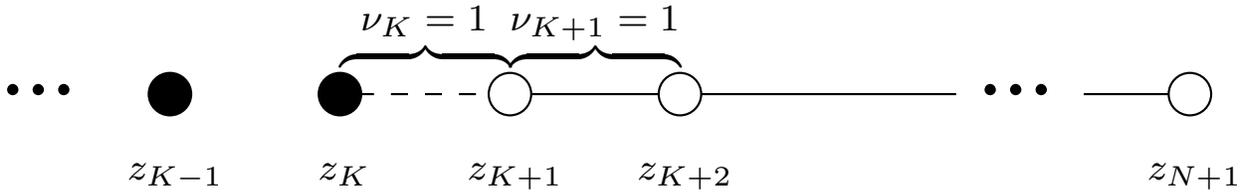}
\caption{\label{qcchain}One end of the quasicontinuum chain, highlighting the
interface. Filled circles are atomistic representative atoms,
whereas the unfilled circles are continuum representative atoms.}
\end{figure}

For simplicity, we will consider a single
atomistic region symmetrically surrounded by continuum regions
large enough that no atomistic degrees of freedom interact with
the surface.  We let the representative atoms in the range
$j = -K+1,\dots,K$ be
atomistic and in the ranges $j=-N,\dots, -K$ and $K+1,\dots,
N$ be continuum.  Figure~\ref{qcchain} depicts one end
of the quasicontinuum chain.
We note that the atomistic model has surface
effects at the two ends of the chain, but the local quasicontinuum
model does not have surface effects.
Thus, this arrangement of
representative atoms with continuum representative atoms
at the ends of the chain will not give surface effects within the QC
approximation.

We assume that $\nu_j=1$ for $j=-K-1,\dots, K+1.$ This guarantees that
$\nu_j=1$ within the second-neighbor
cutoff radius of any atomistic representative atom
and thus allows $F_{-K+1}^a(\zb)$ and
$F_{K}^a(\zb)$ to be computed without interpolation.
The forces are then given by
\begin{align}
\label{qcf}
F_{j}^{QCF}(\zb) &=
\begin{cases}
\left[\phi'(r_{-N} ) + 2 \phi'(2 r_{-N} ) \right], &j = -N,\\
\left[\phi'(r_{j}) + 2 \phi'(2 r_{j})\right] &\\
\quad - \left[\phi'(r_{j-1}) + 2 \phi'(2 r_{j-1})\right],
& -N+1 \leq j \leq -K,\\
\left[\phi'(r_{j}) + \phi'(r_{j} + r_{j+1})\right] &\\
\quad - \left[\phi'(r_{j-1}) + \phi'(r_{j-1} + r_{j-2})\right],
&-K+1 \leq j \leq K,\\
\left[\phi'(r_{j}) + 2 \phi'(2 r_{j})\right] & \\
\quad - \left[\phi'(r_{j-1}) + 2 \phi'(2 r_{j-1})\right],&K+1 \leq j \leq N,
\end{cases}
\end{align}
where
\begin{equation*}
r_j \define D_ja_0=\frac{(z_{j+1}-z_{j })}{ \nu_j}, \qquad j = -N,\dots,N,
\end{equation*}
is the deformed lattice spacing within the $j$th element.

\subsection{An Energy-Based Quasicontinuum Approximation}

There are many quasicontinuum energies~\cite{mill02,e06,jacobsen04}
that can be used
to precondition the iterative solution of the force-based
quasi-continuum approximation~\eqref{precona}.  We will
give an analysis and numerical experiments for the quasicontinuum
energy described in~\cite{mill02} and denoted
here by $\E^{QCE}$ because it seems to be
the simplest to implement and because it converges sufficiently
rapidly.  Here and in the following QCE will refer specifically to
the energy described in~\cite{mill02}, whereas in the introduction it
represented any possible choice of quasicontinuum energy.

QCE assigns an energy to each
degree of freedom according to the model type (atomistic or
continuum), and the sum of all
such energies gives the total QC energy for the chain. We use
the same distribution of atomistic and continuum representative
atoms as above.  Then the
atomistic representative atoms, located in the range $j = -K+1,\dots,K,$
have energy given by
\begin{equation}
\label{atomSplit}
\E^a_j(\rb) \define
\half\Big[\phi(r_j)
+\phi(r_j + r_{j+1})
+ \phi(r_{j-1})
+ \phi(r_{j-1} + r_{j-2})\Big],
\end{equation}
and the continuum representative atoms, located in the range
$j=-N,\dots,-K$ and
$j=K+1,\dots,N+1,$ have energy~\cite{dobs07}
given by
\begin{equation}
\label{contSplit}
\begin{split}
\E^L_j(\rb) &\define
\half \Big[L_jW(D_j)+L_{j-1}W(D_{j-1})\Big]
\end{split}
\end{equation}
where the energy density $W(D_j)$ is
considered to be zero for $j < -N$ or $j > N.$
The quasicontinuum
energy, $\E^{QCE}(\zb),$ for the chain is given by
\begin{equation}
\label{qceTot}
\E^{QCE}(\zb) = \sum_{j = -N}^{-K} \E_{j }^L(\rb) +
\sum_{j=-K+1}^{K} \E_{j}^a(\rb) +
\sum_{j=K+1}^{N+1} \E_{j }^L(\rb).
\end{equation}
In \eqref{force}, we assign forces according to representative atom type
whereas here we have assigned a partitioned energy according to
representative atom type.

We now mention other QC energies, although they will not be used
in the following.  In~\cite{jacobsen04}, the
quasinonlocal method is proposed to attempt to remove the interface inconsistency
by defining a new QC energy.  For this method, special interface atoms are defined
that behave in a hybrid fashion, interacting atomistically with a neighbor
if that neighbor is atomistic, but using the local approximation to determine
the interaction energy otherwise.  For example, if we denote representative
atoms $j=K$ and $K+1$ to be quasinonlocal, then their energy would be
\begin{equation}
\E^Q_j(\rb) =
\half \big[ \phi(r_j) + \phi(2 r_j) + \phi(r_{j-1}) + \phi(r_{j-1} + r_{j-2})
\big].
\end{equation}
However, this method only gives a consistent quasicontinuum energy for a
limited range of interactions (second-neighbor in one dimension), and further
inconsistencies are introduced when attempting to coarsen the continuum region
in higher dimensions~\cite{e06}.  A more general approach that applies to
longer-range interactions is given in ~\cite{e06}, but this approach to the
development of consistent quasicontinuum energies is also currently restricted
to planar interfaces in higher dimensions.
\section{Convergence of the Iterative Method to Solve the QCF Equations}
\label{sec:iter}
\renewcommand{\rb}{\mathbf r}
 We now give a theorem for the convergence of the
iterative algorithm~\eqref{precona} to solve the QCF equilibrium
equations.
Specifically, we give a domain in which the iteration is
well-defined and a contraction. In the following, this contraction
will be an essential portion of the continuation method that is
applied to solve the final equilibrium problem.

Our result extends the theorem in~\cite{dobs07} by allowing the
removal of the hypotheses on the external force, $\tilde{\fb}
\define(\tilde{f}_{-M},\dots,\tilde{f}_M),$ by
utilizing mixed boundary conditions in the problem
analyzed in this paper rather than the free boundary conditions
analyzed in~\cite{dobs07} (the different assumptions lead to
different constants in the inequalities).
 In this
section, the dependence on $s$ of both the solution, $\rb,$ and
external force, $\tilde{\fb},$ is again suppressed.

Since the QCF forces~\eqref{qcf} and the QCE
energy~\eqref{qceTot} depend only on the interatomic spacings,
$\{r_j\}_{j=-N}^N,$ the analysis of the iteration is simplified
by formulating the problem
in terms of forces on the lattice spacing, $\{r_j\}_{j=-N}^N,$ rather than
on representative atom positions, $\{z_j\}_{j=-N}^{N+1}.$
We note that since $z_{N+1}=\hat y_{M+1}(s)$ is fixed, there is a
one-to-one mapping $\zb \leftrightarrow \mathbf r.$ For the
energy-based quasicontinuum approximation, we define
$\psi^{QCE}_j(\rb)$ to be the force conjugate to the
representative atom
spacing $z_{j+1}-z_j=\nu_j r_j,$ namely $\psi^{QCE}_j(\rb)
\define -\nu_j^{-1}\pd{\E^{QCE}}{r_j}(\zb).$
This conjugate force satisfies
\begin{equation*}
-\psi^{QCE}_j(\rb) =
\begin{cases}
  \phi'(r_j) + 2 \phi'(2r_j), &  -N \leq j \leq -K-2,\\
\phi'(r_j) + 2 \phi'(2r_j) + \frac{1}{2} \phi'(r_j +r_{j+1}),
& j = -K-1, \\
\phi'(r_j) + \frac{1}{2} \phi'(r_j+r_{j-1})  & \\
    \qquad + \frac{1}{2} \phi'(r_j+r_{j+1}) + \phi'(2r_j), & j = -K,\\
\phi'(r_j) + \frac{1}{2} \phi'(r_j+r_{j-1}) + \phi'(r_j+r_{j+1}),  &j=-K+1,\\
\phi'(r_j) + \phi'(r_j+r_{j-1}) + \phi'(r_j+r_{j+1}),
       & -K + 2 \leq j \leq K - 2, \\
\qquad \qquad \vdots &
\end{cases}
\end{equation*}
We have from the chain rule that~\cite{dobs07}
\begin{equation}\label{for}
\begin{split}
F^{QCE}_j(\zb)&\define -\pd{\E^{QCE}}{z_j}(\zb)\\
&=-\psi^{QCE}_{j}(\rb)  +
\psi^{QCE}_{j-1}(\rb),\qquad j=-N,\dots,N,
\end{split}
\end{equation}
where $\psi^{QCE}_{-N-1}(\rb)\define 0,$ so it follows by
summing~\eqref{for} that
\begin{equation*}
\psi^{QCE}_j(\rb) = - \sum_{i=-N}^j F^{QCE}_{i}(\zb),\qquad
j=-N,\dots,N.
\end{equation*}
If we define an analogous quantity $\psi^{QCF}_j(\rb)$ by setting
\begin{equation}\label{sumstress}
\psi^{QCF}_j(\rb) \define - \sum_{i=-N}^j F^{QCF}_{i}(\zb),\qquad
j=-N,\dots,N,
\end{equation}
then we have that
\begin{equation*}
- \psi^{QCF}_j(\rb) =
\begin{cases}
\phi'(r_j) + 2 \phi'(2 r_j), & -N \leq j \leq -K ,\\
\phi'(r_j) + \phi'(r_j + r_{j-1}) + \phi'(r_j + r_{j+1})
  + I_{-K},& -K+1 \leq j \leq K,\\
\phi'(r_j) + 2 \phi'(2 r_j) + I_{-K} - I_{K},& K+1 \leq j \leq N,
\end{cases}
\end{equation*}
where
$I_j = 2\phi'(2r_{j}) - \phi'(r_{j} + r_{j-1}) - \phi'(r_{j} + r_{j+1}).$
The external force is likewise made conjugate to the representative atom
spacing by summing
\begin{equation}\label{symm}
\Phi_j = -\sum_{i=-N}^j f_i,\qquad j=-N,\dots,N.
\end{equation}
It follows from~\eqref{sumstress} and~\eqref{symm} that
 a configuration $\zb$ is a solution to
\begin{equation}
\label{equil}
F_{j}^{QCF}(\zb)+f_{j}=0,\qquad j=-N,\dots,N,
\end{equation}
if and only if the corresponding $\rb$ is a solution to
\begin{equation*}
\psi^{QCF}_j(\rb) + \Phi_j = 0,\qquad j=-N,\dots,N.
\end{equation*}

We will iteratively solve the equilibrium equations \eqref{equil} by
using $F^{QCE}$ as a preconditioner for $F^{QCF}.$
The convergence theorem we prove later in this section means
that QCE is quite close to QCF, and the iterative equations
converge rapidly.  Thus, if we use standard
energy minimization algorithms to solve the QCE equations at each
iterative step and utilize the fast convergence of the QCE
solution to QCF, then we get an efficient solution method for the QCF
equations with its inherent advantages of consistency and
simplicity.   The iterative equations are
\begin{equation}
\label{psiitera}
\psi_j^{QCE}(\rb^{p+1}) + \psi_j^{G}({\rb^{p}})
+ \Phi_j = 0, \qquad j=-N,\dots,N,
\end{equation}
where the correction force is
\begin{equation}
\label{gconj}
\psi^G_j(\rb)\define\psi^{QCF}_j(\rb)-\psi^{QCE}_j(\rb).
\end{equation}

\subsection{Assumptions on the Atomistic Potential,  $\phi(r)$}
Before stating our result about the convergence of the
iteration~\eqref{psiitera},
we make explicit the assumptions on the potential, $\phi.$
A prototypical function fitting these assumptions is
the Lennard-Jones potential,
\begin{equation}
\label{lj}
\phi(r) = \frac{1}{r^{12}} - \frac{2}{r^6}.
\end{equation}
We recall that the energy density corresponding
to $\phi(r)$ for the second-neighbor energy~\eqref{atomTot} is given
by
$W(D)=a_0^{-1}\left(\phi(Da_0)+\phi(2Da_0)\right)$ where
$a_0$ is the equilibrium bond length of a uniform chain,
that is,
it is the minimum of $\phi(r)+\phi(2r).$

\begin{figure}
\includegraphics[height=1.75in,width=\textwidth]{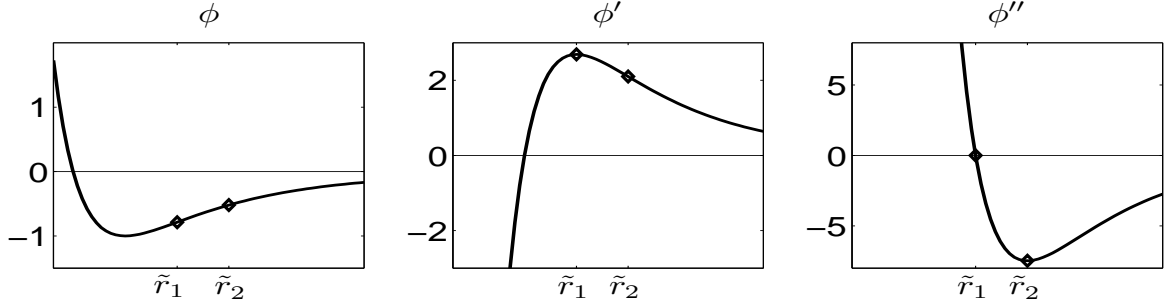}
\caption{\label{phiplot}The Lennard-Jones potential \eqref{lj} demonstrates the
prototypical behavior of $\phi(r)$ and its derivatives.}
\end{figure}

\begin{figure}
\includegraphics[height=1.75in,width=\textwidth]{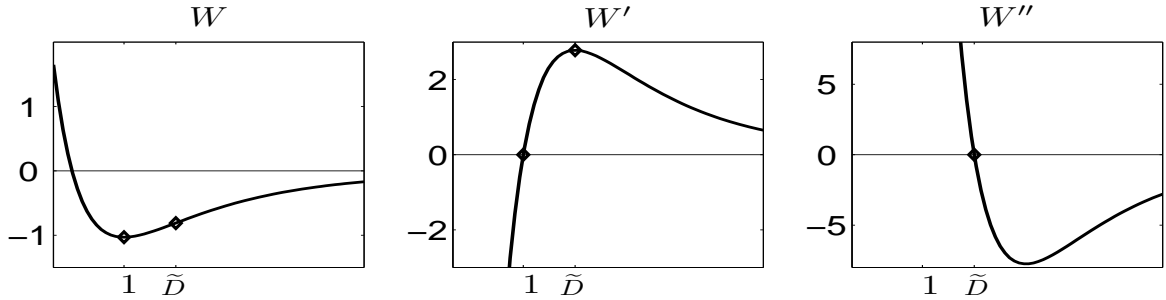}
\caption{\label{hatphiplot}The energy density,
 $W(D),$ corresponding to the Lennard-Jones potential
 \eqref{lj} and its derivatives.}
\end{figure}

We will assume that
$\phi(r) \in C^3\left((0,\infty\right))$
and that it satisfies the following
properties that are illustrated in the Lennard-Jones \eqref{lj}
case in
Figures \ref{phiplot} and \ref{hatphiplot}.  There exist
$\tilde{r}_1,$ $\tilde{r}_2,$ and ${\widetilde{D}}$ such that
\begin{align*}
&\phi''(r) > 0 \text{ for } 0 < r < \tilde{r}_1 \text{ and }
\phi''(r) < 0 \text{ for } r > \tilde{r}_1,\\
&\phi'''(r)<0 \text{ for } 0 < r < \tilde{r}_2
\text{ and  } \phi'''(r)>0 \text{ for } r > \tilde{r}_2, \\
 &W'(D)<0 \text{ for }  0 < D < 1
\text{ and } W'(D)>0\text{ for } D>1,\\
&W''(D)>0 \text{ for }  0 < D<{\widetilde{D}}
\text{ and } W''(D)<0 \text{ for } D>{\widetilde{D}},\\
&0 < a_0 < \tilde{r}_1 < \tilde{r}_2 < 2 a_0, \\
&1 < {\widetilde{D}} .
\end{align*}
We note that ${\widetilde{D}}$
is the deformation gradient of a uniform chain at the load limit.

The following theorem gives sufficient conditions on the existence
of a region $\mathbf r=(r_{-N},\dots,r_N)\in\Omega = (r_L,r_U)^{2N+1}$
in which the iteration \eqref{psiitera} is well-defined and a contraction.
We see that under these conditions
QCE is an efficient preconditioner for
the force-based equations, giving a contraction mapping for the
iteration.
The idea is that these quasicontinuum approximations are
quite close, so that the solution of QCE gives a good approximation to
the solution of QCF.

\begin{theorem}
\label{thm:contract}
For a given conjugate external force, $(\Phi_{-N},\dots,\Phi_{N})$
suppose that there exist $r_L$ and $r_U$ such that
\begin{align}
&\frac{\hat{r}_2}{2} <  r_L < r_U,\label{hat} \\
&\label{rrange} \phi''(r_U) + 21 \phi''(2 r_L) > 0,\\
\phi'(r_L) + 6 \phi'(2 r_L) - 4 \phi'(2 r_U)& < \Phi_j
< \phi'(r_U)
     + 6 \phi'(2 r_U) - 4 \phi'(2 r_L),
     \label{phirange}
\end{align}
for $j = -N,\dots,N.$
Then for every $\mathbf r^{p} \in \Omega := (r_L,r_U)^{2N+1}$
there is a unique $\mathbf r^{p+1} \in \Omega$
such that
\begin{equation}\label{psiiter}
\psi^{QCE}_j(\mathbf r^{p+1}) + \psi^G_j(\mathbf r^{p}) + \Phi_j = 0,
\qquad j=-N,\dots,N.
\end{equation}
We also have that the
induced mapping is a contraction: if $\mathbf r^p \rightarrow \mathbf r^{p+1}$
 and $\mathbf s^p \rightarrow \mathbf s^{p+1},$ then
\begin{equation*} 
\maxnorm{\mathbf r^{p+1} - \sbf^{p+1}}\le \frac{  16 |\phi''(2 r_L)|}
{\phi''(r_U) -  5 |\phi''(2 r_L)|}
\maxnorm{\mathbf r^p - \sbf^p},
\end{equation*}
where we have from \eqref{rrange} that
\begin{equation*}
\frac{  16 |\phi''(2 r_L)|}{\phi''(r_U) -  5 |\phi''(2 r_L)|} < 1.
\end{equation*}
\end{theorem}

We start by remarking on the theorem's assumptions.
 The second inequality in~\eqref{hat}
states that $r_L$ is acting as a lower bound on $\min_j r_j$ and
$r_U$ as an upper bound.  The first inequality is chosen for convenience so
that $\phi''(2 r)$ is monotone.  (We note that for Lennard-Jones and similar
potentials, it is physically a very reasonable assumption
due to the stiffness of the interactions.)
Condition~\eqref{rrange} ensures diagonal
dominance of the Jacobian matrix for $F^{QCE}$ and ensures that
the contraction estimate is less than 1.  Finally, the condition on
$\Phi_j$ in~\eqref{phirange} restricts the external forces sufficiently to
allow a simple degree theory argument to prove existence of solutions.

It is possible to choose a fairly large range for $r$ when the external forces
are far from the load limit of the chain.  However, as $\max_j \Phi_j$
approaches the load limit, $r$ must approach the tensile limit which makes the
estimates much more sensitive.  This reduces the size of
$(r_L,r_U).$
The hypotheses of Theorem~\ref{thm:contract} guarantee that
the iteration~\eqref{psiiter} converges to the QCF approximation
of a stable atomistic solution.

This theorem is a modification of Theorem 5.1 in~\cite{dobs07}.  The
proof there models the QCE equations as a perturbation of the fully local
quasicontinuum energy, $\E^L.$
The proof here follows by modifying the original proof
to handle the new terms that arise from removing the assumption
of symmetry of $\Phi_j.$  These new terms can be estimated by
the techniques used to estimate similar terms
analyzed in~\cite{dobs07}.

\section{The Deformation of a Lennard-Jones Chain under Tension}
\label{tension}
\renewcommand{\rb}{\mathbf r}
In this section, we consider the application of
Theorem~\ref{thm:contract} to a chain modeled by the
Lennard-Jones potential~\eqref{lj}.  The deformation of the
fixed end, $\hat y_{M+1}(s),$ can be set arbitrarily since the
dependence on $\hat y_{M+1}(s)$ is given by a uniform
translation.  To obtain a uniform tension, we model the
external force for the fully atomistic chain by
$\tilde{f}_{-M}=-\Phi$ and $\tilde{f}_j=0$ for $j=-M+1,\dots,M.$ It then follows
from \eqref{symm} that the conjugate external force for the
QC approximation is given by $\Phi_j = \Phi$ for all
$j=-N,\dots,N.$

There are uniform solutions to the
QCF equations up to the load limit $\Phi_{max}=2.7810,$ that is,
if $\phi'(r) + 2 \phi'(2 r) = \Phi,$ then
$\psib^{QCF}(r \eb) = \Phi \eb,$  where
$\eb = (1, 1, \dots, 1) \in\Real^{2N+1}$ and $\psib^{QCF}(\rb)
\define \left(\psi^{QCF}_{-N}(\rb), \dots, \psi^{QCF}_N(\rb) \right).$
We apply an external force very close
to the load limit to get an example where continuation is
necessary to ensure that the preconditioned equations converge.
Define the loading path $\phib(s) =  2.76s \, \eb.$
Then solutions $\rb(s)$ to the QCF equations satisfy
$\rb(s) = r(s) \eb$ where
\begin{equation} \label{rs}
\phi'(r(s)) + 2 \phi'( 2 r(s)) = 2.76 s, \qquad s \in [0, 1].
\end{equation}

\begin{figure}
\includegraphics[width=.49\textwidth]{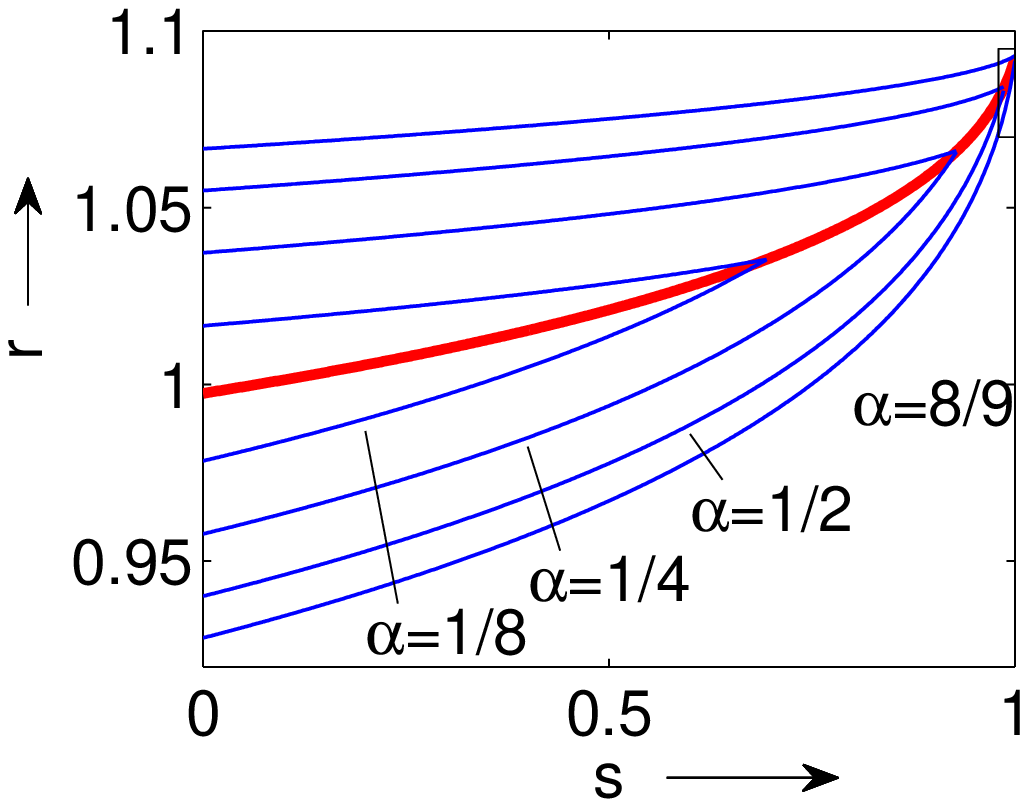}
\includegraphics[width=.49\textwidth]{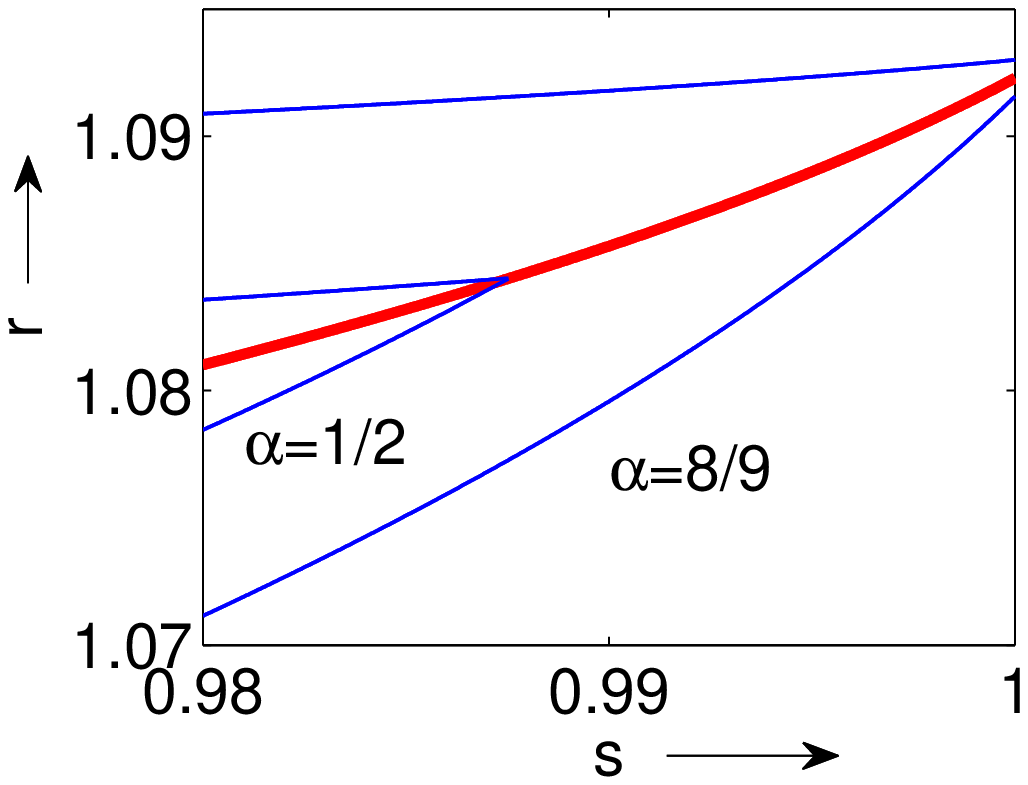}\\
\centerline{\hspace{.2in} (a) \hspace{2in} (b)}
\caption{(a) Loading response, $r(s),$ for the Lennard-Jones chain
is surrounded by contraction regions, $r(s) \pm \delta(s),$
corresponding to $\alpha=\frac{1}{8}, \frac{1}{4}, \frac{1}{2},
\frac{8}{9}$.  The contraction constant $\alpha$ increases with
distance from $r(s).$  (b) Detail shows the contraction region
in a neighborhood of $s=1.$}
\label{alphas}
\end{figure}

For any $s\in[0,1],$ we can apply Theorem \ref{thm:contract}
to this example by picking $r_L$ and $r_U$ such that
\begin{equation}\label{alpha}
\frac{ 16 | \phi''(2 r_L) |}{\phi''(r_U) - 5 | \phi''(2r_L) |} =
\alpha < 1,
\end{equation}
to conclude that the iterative
equation using the QCE preconditioner is a contraction mapping
with contraction constant $\alpha,$ provided that
\eqref{hat}-\eqref{phirange}
holds. We find $r_L$ and $r_U$ symmetrically positioned about $r(s)$
by substituting $r_L(s) = r(s) - \delta(s)$ and $r_U(s) = r(s) + \delta(s)$
into~\eqref{alpha} with $r(s)$ given by the solution to \eqref{rs} to
obtain
\begin{equation}
\label{delta}
\phi''(r(s) + \delta(s)) + (5+16/\alpha)
\phi''(2(r(s)-\delta(s))) = 0.
\end{equation}
It can be checked that \eqref{hat}-\eqref{phirange}
are satisfied with $\Phi_j=\Phi=2.76 s$ for $j=-N,\dots,N.$
Therefore, for any initial guess
 $\rb^0 \in [r(s)-\delta(s), r(s)+\delta(s)]^{2N+1},$
the iteration step \eqref{psiitera} is a contraction
mapping for all $n$ with
contraction rate $\alpha$ and limit point $\rb(s) = r(s) \eb.$
Figure \ref{alphas} depicts the
solution $\rb(s)$ along with four contraction intervals that correspond
to $\alpha = \frac{1}{8}, \frac{1}{4}, \frac{1}{2}, \mathrm{and} \
\frac{8}{9}.$
For every $\alpha\leq 1$ the corresponding
$\delta(s)$ is decreasing.
In Section \ref{sec:contend} we consider the contraction region
corresponding to $\alpha=\frac{8}{9}$ since this contraction
region terminates just beyond our maximum applied load, $\delta(1.001)=0.$

\renewcommand{\rb}{\mathbf r}
\subsection{Fracture at the Interface}
To demonstrate the need for continuation methods
for the above example, we describe the
performance of a modified version of \verb+qc1d+, a code by
Ellad Tadmor for solving the QCF equations using QCE as a
preconditioner with a nonlinear conjugate gradient method for
solving the inner iteration. We attempt to directly solve
\eqref{inner} starting from the energy-minimizing lattice spacing and
using only a single loading step, $Q=1,$ which gives the following
minimization problem.
We have
\begin{equation}\label{onestep}
\mathbf{r}_1^{1}=\argmin_{\rb} \left[\psib^{QCE}(\rb)
-\left(2.76+
\psib^{G}(\mathbf {a_0})\right)\cdot\mathbf {R}(\rb)\right],
\end{equation}
where $\mathbf{R}(\rb) = (\nu_{-N}
r_{-N},\dots,\nu_{N} r_N)$ denotes the representative atom
spacing.

We consider an uncoarsened QC chain, with
\begin{equation*}
M = N = 7,
\end{equation*}
 undergoing
external loading as described in Section~\ref{tension}.  The
chain is partitioned with
\begin{equation*}
K = 3
\end{equation*}
which means that there are six atomistic representative atoms
surrounded symmetrically by two groups of five continuum
representative atoms.  The QCF solution  $\rb(s)$ given
by~\eqref{rs} and the contraction parameters $\delta(s)$ and
$\alpha$ given by~\eqref{delta} do not depend on the size of the
chain; however, the QCE preconditioner solution will depend on
the size and composition of the chain because it has a
non-uniform solution due to the atomistic to continuum
interface.  Because the exact solution is a uniformly deformed
chain, our problem is unchanged by any coarsening of the continuum region.
While this does not illustrate the power of QC approximations to
reduce computational complexity, it provides a simple case in
which to analyze loading up to a singular solution, in this case
fracture.

\begin{figure}
\includegraphics[width=\textwidth]{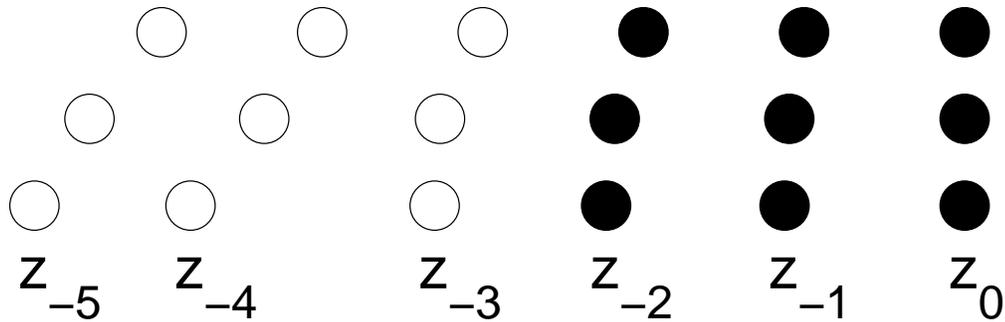}
\caption{\label{fracture}
A close-up of the atomistic to continuum interface
showing fracture that occurs when continuation is not employed.
The three layers represent three steps of a single conjugate
gradient minimization for the iterative
equations~\eqref{onestep}, where the position of $z_0$ has been
normalized to align the chains.  The upper layer is the
initial, undeformed state $\rb = \mathbf{a}_0$.
The middle layer shows a subsequent iteration where
the chain is nearly uniformly deformed and close to the QCF
solution.  The lower layer shows later iteration where a clear
separation of atom pairs occurs.  None of the states shown is
a solution to the minimization problem, and the numerical
algorithm eventually terminates without finding a minimum.}
\end{figure}

The chain fractures in the atomistic to continuum interface.
In the interface, QCE behaves like a continuum material with
varying stiffness which is why it fails to be a consistent scheme.
The corrections $\psib^G(\mathbf{a_0})$
act to counterbalance
this effect by compressing the high stiffness regions, and adding
tension to the low stiffness regions.
Fracture occurs due to the fact that the corrections
\eqref{gconj} applied in the atomistic to continuum transition
are a model correction at the equilibrium bond length, $a_0,$
but are much too strong at the stretched configuration.
The overcorrections add to the very large external force and exceed the load
limit for the QCE chain (see Figure~\ref{fracture}).
The above estimates show that the
continuation method described in Section~\ref{sec:cont} provides a convergent
method for computing the deformation of the chain at the load $\Phi=2.76$ with
the \verb+qc1d+ code.
\section{Solution of the QCF Equations by Continuation}
\label{sec:cont}

In this section and the following two, we give an analysis of the solution
of the QCF equilibrium equations by continuation.  We will
present our results in a general setting that focuses
on the contraction property of the preconditioned equations.  Because we
only use the abstract contraction property, the continuation analysis given
here will apply to higher dimensional QC approximations provided one has
a contraction result similar to Theorem~\ref{thm:contract}.
Given $G \in C^1(\Real^{n+1}; \Real^n),$ our goal will be to
approximate a curve of solutions
$\rb(s)\in C^1([0,1]; \Real^n)$ to
\begin{equation}
\label{eq:me}
G(\rb(s), s) = 0 \qquad \text{ for  } s \in [0,1],
\end{equation}
where
\begin{equation*}
\begin{split}
\det \nr G(\rb(s),s) \neq 0
  &\qquad \text{ for } s \in [0, 1].
\end{split}
\end{equation*}
We will later apply this theory to QCF by considering the solution
of the equations
\[
G(\rb(s),s) \define \psi^{QCF}(\rb(s),s) + \Phi(s) = 0
\qquad \text{ for all } s \in [0,1].
\]

Let $k(s)$ be a bound on $\rb'(s),$ that is,
$\maxnorm{\rb'(s)} \leq k(s),$ which gives
\[
\maxnorm{\rb(t) - \rb(s)} \leq \int^{t}_{s} k(\tau) \, \ud \tau\qquad\text{for
all } t > s.
\]
 We assume further that for each
$s \in [0, 1]$ there is an iterative solver $T_{s}: \Real^n
\rightarrow \Real^n$ that is locally a contraction mapping with
fixed point $\rb(s).$  That is, there is an $\alpha < 1$ such that
for every $s \in [0,1],$ there is a radius $\delta(s)$ with the
property
\begin{equation*}
\maxnorm{\rb(s) - \mathbf p}, \maxnorm{\rb(s)
- \mathbf q} \leq \delta(s) \Rightarrow \maxnorm{T_{s} \mathbf p -
T_{s} \mathbf q} \leq \alpha \maxnorm{\mathbf p - \mathbf q}.
\end{equation*}
We saw in Section~\ref{tension} that such a radius $\delta(s)$
can be obtained for $T_s$ given by the iterative method~\eqref{psiiter} if
the hypotheses of Theorem~\ref{thm:contract} are satisfied.

Let $0 = s_0 \leq s_1 \leq \dots \leq s_Q = 1$ be a sequence of load steps
 where at each point $s_q$ we wish
to compute $\rb_q,$ an approximation to $\rb(s_q).$ Beginning from
an initial guess $\rb^0_q,$ the iterative solver $T_{s_q}$ is
applied to \eqref{eq:me}, keeping $s_q$ fixed.  This generates a
sequence of approximations $\rb^p_q = T_{s_q}^p \rb^0_q$ for $p = 1,
\dots, P_q,$ where $T_{s_q}^p$ denotes $p$ compositions
of the operator
$T_{s_q}$ and $P_q$ denotes the number of iterations at step $m.$
  We
then let $\rb_q = \rb^{P_q}_q.$ The choice of initial guess
$\rb^0_q$ is typically made using polynomial extrapolation, and
here we choose $r^0_q = r_{q-1},$ which is zeroth-order
extrapolation.

We will now give an analysis of how to choose
the load steps $0 = s_0 \leq s_1 \leq \dots \leq s_Q = 1$ and
the corresponding number of iterations $P_1,\dots,P_Q$ to efficiently
approximate $\rb(s)$ with respect to two different goals.
We first consider the efficient approximation of
$\rb(s)$ in the maximum norm
for all $s\in[0,1],$ and we then consider the
efficient approximation of the end point $\rb(1).$
We note that our analysis only gives
an upper bound for the amount of work needed to
compute an approximation of our chosen goal
to a specified tolerance
since we use a uniform estimate for the rate
of convergence $\alpha$
rather than the decreasing $\alpha$ as we
converge to the solution that we can obtain from
Theorem~\ref{thm:contract} and displayed in
Figure~\ref{alphas}.

\section{Efficient Approximation of the Solution Path
in the Maximum Norm}
\label{sec:uniform}
For simplicity, we will first consider a
uniform region of contraction radius $\delta(s)=\delta,$
a uniform bound $k(s)=k,$
a uniform step size $h=1/Q= s_q-s_{q-1},$
and a uniform number of iterations
at each step $P_q=P.$  We will denote the continuous,
piecewise linear interpolant of
$\rb(s_q)\in\Real^n,
\ q=0,\dots,Q,$ by $\I\rb(s);$
and we will denote the continuous,
piecewise linear interpolant of
$\rb_q\in\Real^n,
\ q=0,\dots,Q,$ by $\tilde \rb(s).$
We will determine an efficient
choice of $h$ and $P$ to guarantee that
\begin{equation}\label{guar}
\max_{s\in[0,1]}\maxnorm{\rb(s)-\tilde \rb(s)}\le2\epsilon,
\end{equation}
where we assume for convenience that $2 \epsilon \le \delta.$

We will assume that $\rb(s)\in C^2([0,1]; \Real^n),$
so there exists a constant $k_2\ge 0$ such
that
\begin{equation*}
\max_{s\in[0,1]}\maxnorm{\rb(s)-\I\rb(s)}\le k_2h^2.
\end{equation*}
We can then ensure that
\begin{equation*}
\max_{s\in[0,1]}\maxnorm{\rb(s)-\I\rb(s)}\le \epsilon
\end{equation*}
by choosing $h\le \sqrt{\epsilon/k_2}.$  We can thus satisfy
\eqref{guar} by guaranteeing that
\begin{equation}\label{guar2}
\max_{q=0,\dots,Q}\maxnorm{\rb(s_q)-\rb_q}\le\epsilon.
\end{equation}

Now if $\maxnorm{\rb(s_{q-1})-\rb_{q-1}}\le\epsilon$, then
\begin{align*}
\maxnorm{\rb(s_q)-\rb^0_{q}}
&\le \maxnorm{\rb(s_q)-\rb(s_{q-1})}+\maxnorm{\rb(s_{q-1})-\rb_{q-1}}\\
&\le kh+\epsilon.
\end{align*}
We choose $0<h\le \frac{\delta-\epsilon}{k}$ so that $\rb^0_q$ is in the
region of contraction
\[
\maxnorm{\rb(s_q)-\rb^0_q}\le kh+\epsilon\le\delta.
\]
We then choose $P$ to achieve the desired error
\begin{align*}
\maxnorm{\rb(s_q)-\rb^P_q}&\le\alpha^P\maxnorm{\rb(s_q)-\rb^0_q}\\
&\le\alpha^P(kh+\epsilon)\le\epsilon.
\end{align*}
We can thus guarantee that
$\maxnorm{\rb(s_q)-\rb^P_q}\le\epsilon$ by
doing $P$ iterations where
\[
 P(h)=\frac   {\ln \left(\frac{\epsilon}{\epsilon+kh}\right)}  {\ln \alpha}.
\]

The computational work to obtain \eqref{guar2} can then
 be bounded by
\[
\W(h)=\frac {P(h)}{h}=\frac   {\ln
\left(\frac{\epsilon}{\epsilon+kh}\right)}{h\ln
\alpha} \quad \text{for }0\le h\le
\frac{\delta-\epsilon}{k}.
\]
We have by the Mean Value Theorem that
\begin{equation*}
\frac{d \W(h)}{dh}= \frac{k}{h\ln\alpha} \left[
\frac{\ln(\epsilon+kh)-\ln\epsilon}{kh}
-\frac1{(\epsilon+kh)} \right]<0
\end{equation*}
for $ 0<h< \frac{\delta-\epsilon}{k}.$

We can finally obtain \eqref{guar} by choosing
\begin{align*}
h_{opt}&=\min
\left\{\frac{\delta-\epsilon}{k},\sqrt{\epsilon/k_2}\right\} , \\
P(h_{opt})&=\frac   {\min \left\{ \ln \frac{\epsilon}{\delta},
\ln \frac{\sqrt k_2\epsilon}{\sqrt k_2\epsilon+{k}\sqrt\epsilon} \right\} }
 {\ln \alpha}
\to \infty\quad\text{as}\quad \epsilon\to 0.
\end{align*}
As $\epsilon$ goes to zero, the second criterion becomes active
so that the step size is determined by the interpolation
estimates rather than the size of the contraction region.  The
number of steps grows as $\frac{1}{\sqrt{\epsilon}},$ and the
number of iterations per step grows like like $\frac{\ln \epsilon}{2 \ln
\alpha}.$
\section{Efficient Approximation of the Solution
at the Final State}
\label{sec:contend}

In this section, our goal will be to compute
$\rb_Q$ satisfying the error tolerance
\begin{equation}
\label{error}
\maxnorm{\rb(1) - \rb_Q} \leq \epsilon
\end{equation}
while minimizing the computational effort
\begin{equation*}
\W(\{P_q\},\{s_q\}) \define\widehat \W \sum_{q=1}^Q P_q,
\end{equation*}
subject to the constraints on $\{P_q\}$ and $\{s_q\}$ given below,
where $\widehat \W>0$ is the work per iterative step which we scale to
$\widehat \W=1.$  We note that the preceding assumes that applying the
iterative solver is the most computationally expensive operation
and all iterations are equally expensive. We first formulate the
optimization problem with constraints, and we then simplify the
problem by observing that some of the inequality constraints can
be replaced by equality constraints.  In this section, we now
consider a continuous, decreasing contraction radius $\delta(s)$
and a continuous, positive bound $k(s).$  The load steps taken to
achieve the error goal will not be uniformly spaced which will
take advantage of the large initial contraction region and the
fact that low error is only desired at the endpoint, $s=1.$

We define the error at $s_{q}$ by $e_{q} =
\maxnorm{\rb(s_{q})-\rb_{q}}$ for all $q = 0,\dots,Q.$ Then a
bound on the error in the initial guess $\rb^0_{q}$ for $q=1,\dots,Q$
can be given by
\begin{equation*}
\begin{split}
\maxnorm{\rb(s_{q})-\rb^0_{q}}
        &\le \maxnorm{\rb(s_{q})-\rb(s_{q-1})}
            + \maxnorm{\rb(s_{q-1})-\rb_{q-1}} \\
        &\le \K(s_{q}) - \K(s_{q-1}) + e_{q-1} ,
\end{split}
\end{equation*}
where $\K(s) = \int_0^s k(\tau) \, \ud \tau.$ If
$\K(s_q) - \K(s_{q-1}) + e_{q-1} \leq \delta(s_{q}),$
the mapping $T_{s_{q}}$ is a
contraction and the error satisfies the bound
\begin{equation*}
e_{q} = \maxnorm{\rb(s_{q}) - \rb^{P_{q}}_{q}} \\
    \leq \alpha^{P_{q}}\maxnorm{\rb(s_{q})-\rb^0_{q}}
\end{equation*}
since $\rb(s_{q})$ is a fixed point of $T_{s_{q}}.$ We assume that
$e_0<\delta(0),$ and we let $\{\gamma_q\}^{Q}_{q=0}$ be the supersolution
for $\{e_{q}\}^{Q}_{q=0}$ defined by the recurrence
\begin{equation*}
\begin{split}
\gamma_q &= \alpha^{P_q} \left(\K(s_{q}) - \K(s_{q-1}) +
\gamma_{q-1}\right),
\qquad \textrm{ for } q = 1,\dots,Q, \\
\gamma_0 &= e_0.
\end{split}
\end{equation*}
In the following, we satisfy the error goal \eqref{error}
 by making sure that the supersolution satisfies $\gamma_Q
\leq \epsilon.$

\begin{figure}
\includegraphics[width=.49\textwidth]{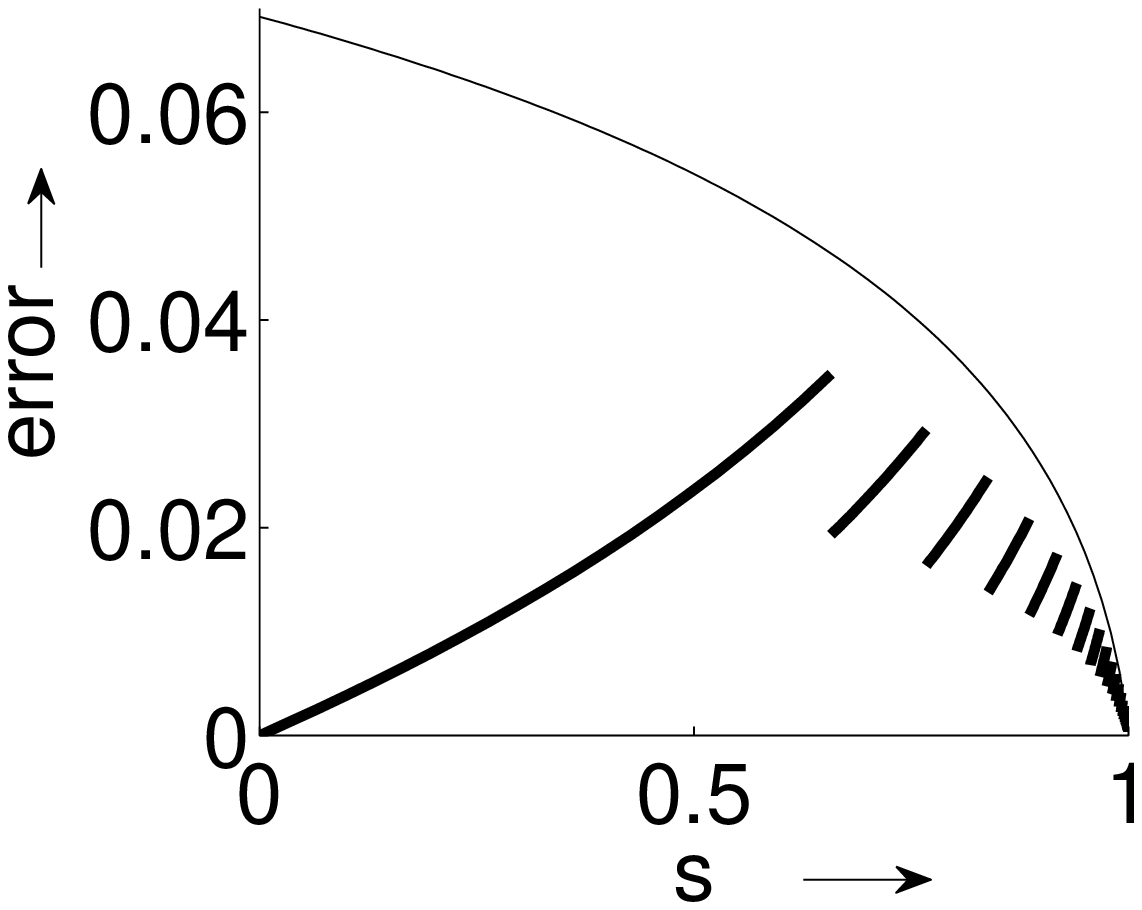}
\includegraphics[width=.49\textwidth]{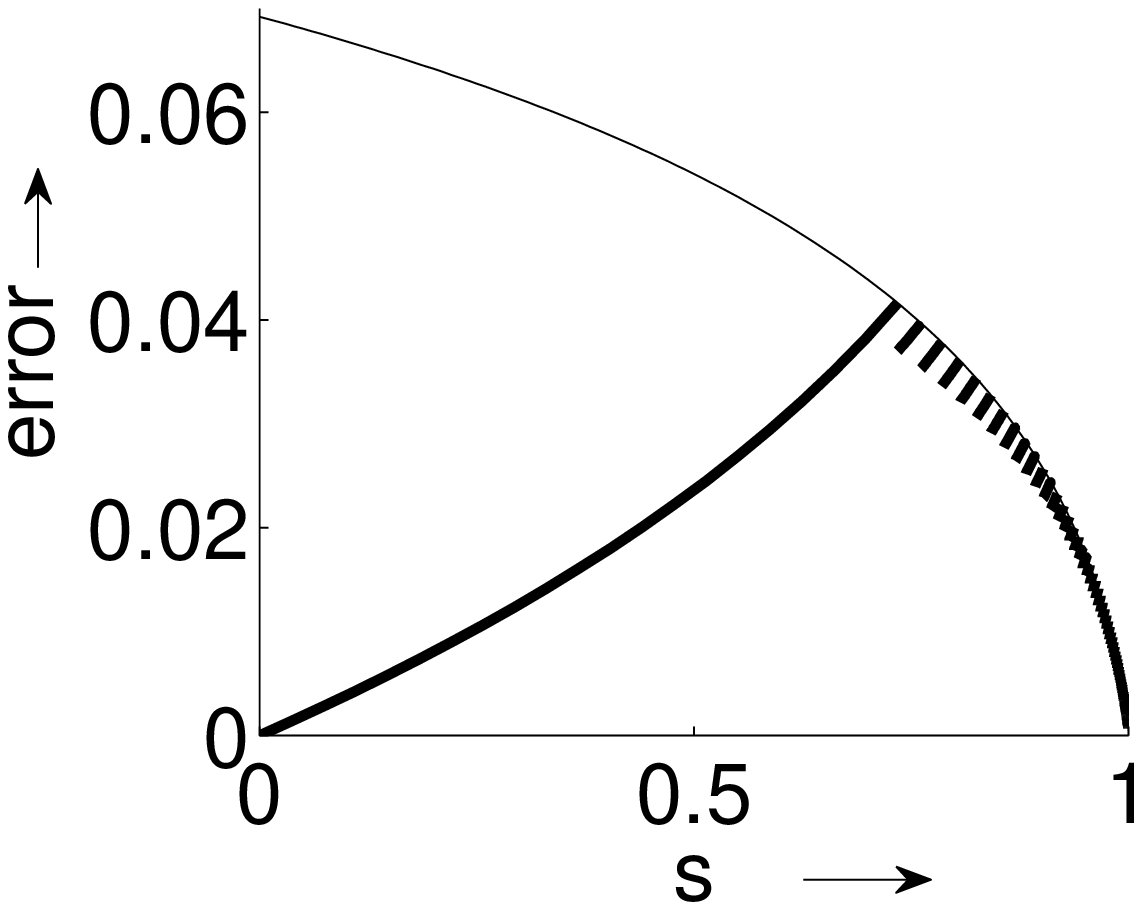}\\
\centerline{\hspace{.2in} (a) \hspace{2in} (b)}
\caption{\label{loading}
Error $e(s) = \maxnorm{\rb(s) - \hat{\rb}(s)}$ for the deformation
problem given in Section~\ref{tension},
where
$\hat{\rb}(s) = \rb_{q-1},\ s \in [s_{q-1}, s_q).$
Two example
admissible loading paths \eqref{admissible} are displayed in bold.
The contraction radius $\delta(s)$ corresponding to
$\alpha=\frac{8}{9}$ bounds the estimated error curve
for any admissible loading path. (a) Generic path with multiple
iterations per step.  (b) Path with a single iteration per
step and error estimate just less than $\delta(s).$  Path is
optimal solution for Problem \ref{simple}.}
\end{figure}

We now consider the question of how to achieve the error goal
for the supersolution while using the fewest possible applications
of the iterative
solver. We define the set of admissible loading paths that satisfy
the preceding theory by
\begin{equation}
\label{admissible}
\begin{split}
\A = \bigcup^{\infty}_{Q=1} \Big\{ &\left( \{P_q\}^Q_{q=1}, \{s_q\}^Q_{q=0}
\right)
    \subset \Int^{Q-1}_{>0} \times \Int_{\geq 0} \times [0,1]^{Q+1} : \\
 &\qquad 0=s_0 \leq s_1 \leq \dots \leq s_Q = 1,  \\
 &\qquad \K(s_q) - \K(s_{q-1}) + \gamma_{q-1} \leq \delta(s_{q}) \\
  &\qquad \textrm{for all } q = 1,\dots,Q,\text{ and }\gamma_Q \leq \epsilon
 \Big\}.
\end{split}
\end{equation}
Figure \ref{loading} shows the error $\maxnorm{\rb(s) - \hat{\rb}(s)}$
for hypothetical loading paths
\begin{equation*}
\begin{split}
\hat{\rb}(s) = \rb_{q-1} \quad \textrm{for } s \in [s_{q-1}, s_q),
\end{split}
\end{equation*}
using the worst-case error bound
\begin{equation*}
\rb(s_q) - \rb_q = \alpha^{P_q} ( \rb(s_q) - \rb_{q-1}).
\end{equation*}
We will next consider minimizing the work with respect to all
admissible paths.
\begin{problem}
\label{simple}
Given $\epsilon > 0,$ $\K(s) \geq 0$ continuous and increasing,
$\delta(s)>0$ continuous and decreasing,  $0\le \gamma_0
<\delta(0),$ and $0<\alpha<1,$  find
\begin{equation*}
\argmin_{(\{P_q\}, \{s_q\}) \subset \A}  \W(\{P_q\}, \{s_q\})
= \argmin_{(\{P_q\}, \{s_q\}) \subset \A} \ \sum_{q=1}^Q P_q.
\end{equation*}
\end{problem}
We can see that $\A$ is non-empty by considering
paths with sufficiently many small steps so that the error stays
within the contraction domain of the iteration.  Thus, the problem has a
minimizer since
the work for each path is a positive integer.  We denote the
minimum possible work by $\W_{min}.$
The above problem can be analytically solved by the using following
two lemmas which characterize paths that are optimal
in the sense of this problem.
  We first observe that it is
optimal to only do enough work to stay within the contraction
bounds, that is, the minimum work lies on the boundary
$\K(s_{q}) - \K(s_{q-1}) + \gamma_{q-1} = \delta(s_q).$

\begin{lemma}\label{max}
Let $\LP = (\{P_q\}^Q_{q=1}, \{s_q\}^Q_{q=0}) \in \A$ denote an admissible
loading path.  Then there is $\widehat{\LP} =
(\{\widehat{P}_q\}^Q_{q=1}, \{\widehat{s}_q\}^Q_{q=0}) \in A$  and $J,\,
1 \leq J \leq Q,$ such that $\widehat{P}_q = P_q$ for all
$q = 1,\dots,Q;$
\begin{equation}
\label{maxstep}
\K(\widehat{s}_q) - \K(\widehat{s}_{q-1}) + \widehat{\gamma}_{q-1} =
\delta(\widehat{s}_q)
\end{equation}
for every $q = 1,\dots,J-1;$ and $\widehat{s}_q = 1$ for every $q > J.$
Furthermore, $\widehat{\gamma}_Q \leq \gamma_{Q}$ with equality if
and only if $\LP = \widehat{\LP}.$
\end{lemma}
The full proof is given in the Appendix.  The idea is that since
our goal is to only get accuracy at $s_Q$~\eqref{error}, reducing
error early results in extra total work.  If~\eqref{maxstep}
is not satisfied at some $s_q,$ then we can take a larger step
between $s_{q-1}$ and $s_q$ and smaller
steps later (for $q = 1,\dots,Q-1$) which reduces the
supersolution for the error for all subsequent steps.

Next, we denote the set of all admissible loading paths
with work $\sum^Q_{q=1} P_q \leq m $ by
\begin{equation*}
\A_m = \Big\{ \left( \{P_q\}^Q_{q=1}, \{s_q\}^Q_{q=0} \right) \in
\A : \sum^Q_{q=1} P_q \leq m \Big\} \quad \textrm{ for }
m \geq 0.
\end{equation*}
If the minimum total work is denoted by $\W_{min},$ then
$\A_{\W_{min}}$ is non-empty.  We now show that it is optimal to
only do one iteration per step, increasing the number of steps
as necessary.

\begin{lemma}
\label{splitstep}
There is a path $(\{P_q\}^Q_{q=1}, \{s_q\}^Q_{q=0}) \in
\A_{\W_{min}}$ such that $P_q = 1$ and~\eqref{maxstep} hold for
all
$q=1,\dots,Q-1.$  This uniquely determines $Q.$  Furthermore,
this path has the lowest estimated error, $\gamma_Q,$ of all paths in
$\A_{\W_{min}}.$
\end{lemma}
The full proof is given in the Appendix where it is shown that
any step other than the last with $P_q >1$ can be split to create a new
 admissible path that has the same total work and a lower error.

Combining these two lemmas, we can characterize the optimal path
for solving Problem \ref{simple} by the following algorithm,
where some equations are given in implicit form:
\begin{center}
\begin{alltt}
let \(\gamma\sb{0}=e\sb{0}, s\sb{0}=0, q=0\)
while \(s\sb{q}<1\)
    \(q = q+1\)
    \(\tilde{s}\sb{q} =\) solve\((\K(\tilde{s}\sb{q})-\K(s\sb{q-1})+\gamma\sb{q-1} = \delta(\tilde{s}\sb{q}))\)
    \(s\sb{q} = \min(\tilde{s}\sb{q},1)\)
    \(P\sb{q} = 1\)
    \(\gamma\sb{q} = \alpha\delta(s\sb{q})\)
end
\(P\sb{Q} = \left\lceil \frac{\log\epsilon-\log(\K(1)-\K(s\sb{Q-1})+\gamma\sb{Q-1})}{\log\alpha} \right\rceil\)
\end{alltt}
\end{center}
where $\lceil x \rceil$ is the least integer not less than $x.$

Figure~\ref{loading}b depicts the loading
curve and optimal loading path for our example, where
we directly use
\begin{equation*}
\K(s) \define \int^{s}_0 \maxnorm{ \rb'(\tau)} \, \ud \tau
\end{equation*}
in Problem~\ref{simple}.  We note that the solution depicted in
Figure~\ref{loading}b uses information about the exact
solution, both in the growth estimate $\K$ and in the
computation of contraction regions~\eqref{delta}.
In practice, the results will be applied using estimates to determine
$\delta$ and $\K.$  The lemmas provide the general intuition that
it is efficient to do many steps with a single iteration per step, rather
than fewer steps with many iterations per step.

\appendix
\section{Proofs of Lemma~\ref{max} and
Lemma~\ref{splitstep}}
\label{appendix}
 We present detailed proofs of Lemma~\ref{max} and
Lemma~\ref{splitstep}, which use very similar estimates to show that a
given new path is computationally favorable.

{\em Proof of Lemma~\ref{max}.}
Let $J$ be the smallest integer such that $s_J = 1.$   If $J=1,$
then we are
done; otherwise we show that if
\eqref{maxstep} does not hold for some $q=1,\dots,J-1,$ then we can
adjust $\{s_q\}$ to satisfy \eqref{maxstep} with
a strict decrease in the total error.

Let $j$ be the smallest integer such that
$\K(s_{j}) -\K(s_{j-1}) + \gamma_{j-1} < \delta(s_j).$
If $j < J,$ then our step was too conservative so we
define a new stepping path $\{\widetilde{s}_q\}_{q=0}^{Q}.$
Let $\Delta s > 0$ be chosen so that
$\K(s_{j}+\Delta s) - \K(s_{j-1}) + \gamma_{j-1} = \delta(s_j).$
We let $\Delta s = \min(\Delta s, 1-s_j)$ and
define $\widetilde{s}_q$ by
\begin{equation*}
\widetilde{s}_q = \begin{cases}
s_q,                       & q < j, \\
s_j + \Delta s,            & q = j, \\
\max(s_q, s_j + \Delta s), & q > j.\\
\end{cases}
\end{equation*}
This gives a new loading path
$0 = \tilde{s}_0 \leq \tilde{s}_1 \leq \dots \leq
\tilde{s}_Q = 1.$
By construction, steps $q=1,\dots,j$ satisfy~\eqref{maxstep},
but we must show that all subsequent steps are inside the
contraction region, that is
$\K(\tilde{s}_{q}) -\K(\tilde{s}_{q-1}) +\gamma_{q-1} \leq
\delta(\tilde{s}_q)$ for $q > j.$ If $s_{j+1}> s_j + \Delta s,$ we have
\begin{equation*}
\begin{split}
\widetilde{\gamma}_{j+1}  &=
\alpha^{P_{j+1}} [\K(\tilde{s}_{j+1}) - \K(\tilde{s}_j) +
\widetilde{\gamma}_{j}]\\
&= \alpha^{P_{j+1}} [\K(s_{j+1}) - \K(s_{j} + \Delta s)
+ \alpha^{P_j} \left( \K(s_{j}+\Delta s) - \K(s_{j-1})
+ \gamma_{j-1} \right)] \\
&= \alpha^{P_{j+1}} [\K(s_{j+1}) -\K(s_{j} + \Delta s)
+ \alpha^{P_j} (\K(s_{j}+\Delta s) - \K(s_{j}))  \\
&\qquad + \alpha^{P_j}
\left(\K(s_{j}) - \K(s_{j-1})
+ \gamma_{j-1} \right)] \\
&< \alpha^{P_{j+1}} [\K(s_{j+1}) - \K(s_j) + \gamma_j]\\
&= \gamma_{j+1}.
\end{split}
\end{equation*}
The above shows that the error supersolution $\gamma_{j+1}$ is
reduced and, by consideration of the terms in brackets,
that  the approximation is inside the contraction region at
$s_{j+1}.$  A similar argument holds
for the first non-degenerate step, $s_{j^*} > s_j,$ in the case $s_j+\Delta
s>s_{j+1}.$
Since the remaining path is unchanged, we have $\tilde{\gamma}_q
< \gamma_q$ for all $q = j,\dots,Q.$  We continue this process
until the hypotheses are satisfied. \hfill $\square$

{\em Proof of Lemma~\ref{splitstep}.}
We choose a path in $\A_{\W_{min}}$ of the form given by Lemma~\ref{max}.
Now, if $s_J=1$ for $J<Q,$ we can
combine step $J$ and $J+1$ by letting $\widetilde{P}_{J} =
P_{J+1} + P_J,$ thus we can consider paths
where~\eqref{maxstep} holds for all
 $q = 1, \dots, Q-1.$

Now, we show that if $P_j>1$ for some $j<Q-1,$ then the
error can be reduced by adding a new load step between $s_j$ and
$s_{j+1}.$
Suppose the path satisfies~\eqref{maxstep} for all $q=1,\dots,Q-1$
and $P_j > 1$ for some $j < Q-1.$
We will consider the new path $(\{\tilde{P}_q\}^{Q+1}_{q=0},
\{\tilde{s}_q\}^{Q+1}_{q=1}) \in \A_{\W_{min}}$ given by
\begin{equation*}
\begin{split}
\widetilde{P}_q &= \begin{cases}
P_q       &q < j \\
1     &q = j \\
P_j - 1 &q = j+1 \\
P_{q-1}   &q > j+1
\end{cases}\\
\tilde{s}_q &= \begin{cases}
s_q                & q < j+1 \\
s_j+ \Delta s           & q = j+1 \\
s_{q-1}            & q > j+1,
\end{cases}
\end{split}
\end{equation*}
where $\Delta s$ is chosen such that
\begin{equation*}
\begin{split}
\K(\tilde{s}_{j+1}) &- \K(\tilde{s}_{j}) + \tilde{\gamma}_{j} \\
&= \K(s_{j} + \Delta s) - \K(s_{j}) + \alpha \delta(s_{j})\\
&= \delta(s_{j} + \Delta s) \\
&= \delta(\tilde{s}_{j+1}).
\end{split}
\end{equation*}
The above has a solution, with $0 < \Delta s < s_{j+1}-s_j,$ by the
Intermediate Value Theorem.  The above path clearly has the same
total work as the original, and we now show that it satisfies
the
contraction region constraints in Problem \ref{simple}.
We find that
\begin{equation*}
\begin{split}
\K(\tilde{s}_{j+2}) &- \K(\tilde{s}_{j+1}) + \tilde{\gamma}_{j+1} \\
&= \K(\tilde{s}_{j+2}) - \K(\tilde{s}_{j+1})
   + \alpha^{\widetilde{P}_{j+1}} \delta(\tilde{s}_{j+1}) \\
&= \K(s_{j+1}) - \K(s_{j} + \Delta s) \\
&\qquad + \alpha^{P_{j}-1}
    \left( \K(s_{j} + \Delta s) - \K(s_{j})
    + \alpha \delta(s_{j}) \right)\\
&< \delta(s_{j+1}).
\end{split}
\end{equation*}
Thus, we have lowered the supersolution for the error.  We can
apply Lemma \ref{max} and the above argument until the
hypotheses are satisfied, and at each step the error
supersolution is reduced.
\hfill $\square$

\end{document}